\documentclass[a4paper,reqno]{amsart}
\usepackage{enumerate,amsmath,amsthm,amsfonts}
\usepackage{amsmath}
\usepackage{amsfonts}
\usepackage{amssymb}

\usepackage{epsfig}
\usepackage{graphicx,color}
\usepackage{graphpap,color}
\usepackage{graphicx}
\usepackage{latexsym}
\usepackage{float}
\usepackage{setspace}
\usepackage{lineno}   

\newtheorem{theorem} {Theorem}[section]

.360pk
\font\rt=cmss9.360pk
\font\sd=cmcsc9.360pk

\include {mak}
\begin{document}
~\vspace{-16mm}

\newcommand{\R}{\mathbb{R}}
\newcommand{\Z}{\mathbb{Z}}
\newcommand{\cH}{\mathcal{H}}

\newcommand{\norm}[1]{\left\|#1\right\|}
\newcommand{\hkd}[1]{\left\langle#1\right\rangle}
\newcommand{\krk}[1]{\left\{#1\right\}}
\newcommand{\krb}[1]{\left(#1\right)}
\newcommand{\abs}[1]{\left|#1\right|}
\newcommand{\ol}[1]{\overline{#1}}
\newcommand{\rref}[1]{[\ref{#1}]}

\oddsidemargin 16.5truemm
\evensidemargin 16.5truemm

\thispagestyle{plain}

\vspace{1.5cc}

\begin{center}
{\Large\bf ON GENERALIZED H\"{O}LDER'S INEQUALITY IN WEAK MORREY SPACES}\\ 
\vspace{1.5cc}
{\large\sc ASYRAF WAJIH$^{1}$ and HENDRA GUNAWAN $^{2}$}\\

\vspace{0.3 cm} {\small Faculty of Mathematics and Natural Sciences,\\ Institut Teknologi Bandung,
Bandung 40132, Indonesia\\ $^{1}$asyraf.wajih@hotmail.com, $^{2}$hgunawan@math.itb.ac.id}

\rule{0mm}{6mm}\renewcommand{\thefootnote}{}\footnotetext{
\\
}

\vspace{1.5cc}

\parbox{24cc}{{\Small{\bf Abstract.}
In this note we reprove generalized H\"{o}lder's inequality in weak Morrey spaces.
In particular, we get sharper bounds than those in \cite{gunawan2}. The
bounds are obtained through the relation of weak Morrey spaces and weak
Lebesgue spaces.
}}
\end{center}

\vspace{0.25cc}
\parbox{24cc}{\Small {\it Key words and phrases}: {\it generalized H\"{o}lder's
inequality, weak Lebesgue spaces, weak Morrey spaces.}
}

\vspace{1.5cc}

\section{INTRODUCTION}

H\"{o}lder's inequality is one of the classic inequality which is proved by L.J. Rogers
in 1888 and reproved by O. H\"{o}lder in 1889. Many researchers have obtained new
results related to H\"{o}lder's inequality. Recently, Ifronika \textit{et al.} \cite{gunawan2}
obtained sufficient and necessary conditions for generalized H\"{o}lder's inequality
in generalized Morrey spaces. Furthermore, they also get similar results for weak
Morrey spaces and generalized weak Morrey spaces. In this paper, we present new
bounds for generalized H\"{o}lder's inequality in weak Morrey spaces which are sharper
than those in \cite{gunawan2}.

For $1 \leq p \leq q < \infty$, the \textbf{weak Morrey space} $w\mathcal{M}^p_q(\mathbb{R}^n)$ is
the set of all 
measurable functions $f:\mathbb{R}^n \rightarrow \mathbb{R}$ such that
$$
\|f\|_{w\mathcal{M}^p_q} = \sup\limits_{B = B(a,r), \gamma>0} \left| B \right|^{\frac{1}{q}
- \frac{1}{p}} \gamma \left| \{x \in B: |f(x)| > \gamma \} \right| < \infty.
$$
Here $|B|$ denotes the Lebesgue measure of open ball $B = B(a,r)$ centered at $a$
with radius $r$.
Note that $\|\cdot\|_{w\mathcal{M}^p_q}$ defines a quasinorm on $w\mathcal{M}^p_q(\mathbb{R}^n)$.
Furthermore, if $p = q$, then $w\mathcal{M}^p_q(\mathbb{R}^n) = L^{p, \infty}(\mathbb{R}^n)$.
Thus, $w\mathcal{M}^p_q(\mathbb{R}^n)$ can be considered as generalization of weak Lebesgue
space $L^{p,\infty}(\mathbb{R}^n)$. Consequently, the quasinorm $\|\cdot\|_{w\mathcal{M}^p_q}$
can be rewritten as
$$
\|\cdot\|_{w\mathcal{M}^p_q} = \sup\limits_{B = B(a,r)} \left| B \right|^{\frac{1}{q} -
\frac{1}{p}} \|\cdot\|_{L^{p,\infty}(B)},
$$
where $\|\cdot\|_{L^{p,\infty}(B)}$ is a quasinorm on weak Lebesgue spaces $L^{p,\infty}(B)$.

From \cite{gunawan1}, we have the inclusion relation between weak Morrey spaces
$w\mathcal{M}^{p_1}_q(\mathbb{R}^n)$ and $w\mathcal{M}^{p_2}_q(\mathbb{R}^n)$ for
$1 \leq p_1 \leq p_2 \leq q < \infty$, as stated in the following theorem.

\vspace{.5cc}

\begin{theorem}\label{inklusimorreylemah}\cite{gunawan1}
	If $1 \leq p_1 \leq p_2 \leq q < \infty$, then $w\mathcal{M}^{p_2}_q(\mathbb{R}^n)
\subseteq w\mathcal{M}^{p_1}_q(\mathbb{R}^n)$ with
	$$
\|f\|_{w\mathcal{M}^{p_1}_q} \leq \|f\|_{w\mathcal{M}^{p_2}_q}
$$
	for every $f \in w\mathcal{M}^{p_2}_q(\mathbb{R}^n)$.
\end{theorem}\vspace{0.2mm}

\vspace{1.5cc}

\section{MAIN RESULTS}

Let us first state sufficient and necessary conditions for generalized H\"{o}lder's
inequality in weak Morrey spaces \cite{gunawan2}.

\vspace{.5cc}

\begin{theorem}\label{holdermorreygunawan} \cite{gunawan2}
Let $m \geq 2$. If $1 \leq p \leq q < \infty$ and $1 \leq p_i \leq q_i < \infty$,
$i= 1, 2, \dots, m$, then the following statements are equivalent:

\noindent{\rm (i)} $\sum_{i=1}^m \frac{1}{p_i} \leq \frac{1}{p}$ and $\sum_{i=1}^m \frac{1}{q_i} = \frac{1}{q}$.

\noindent{\rm (ii)} $\biggl\| \prod_{i=1}^{m} f_i \biggr \|_{w\mathcal{M}^{p}_{q}} \leq m \prod_{i=1}^{m}
\| f_i \|_{w\mathcal{M}^{p_i}_{q_i}}$, for every $f_i \in w\mathcal{M}^{p_i}_{q_i}(\mathbb{R}^n)$, $i = 1, 2, \dots, m$.
\end{theorem}\vspace{0.2mm}

\vspace{.5cc}

Now we will refine the statement of Theorem \ref{holdermorreygunawan}, particularly the
part which states that (i) implies (ii). Precisely, we will replace the bound $m$ with a smaller constant.

\vspace{.5cc}

\begin{theorem}\label{holdermorreylemah}
Let $m \geq 2$, $1 \leq p \leq q < \infty$, and $1 \leq p_i \leq q_i < \infty$, $i = 1, 2, \dots, m$.
If $\frac{1}{p^*} = \sum_{i=1}^{m} \frac{1}{p_i} \leq \frac{1}{p}$ and $\sum_{i=1}^{m} \frac{1}{q_i} =
\frac{1}{q}$, then for every $f_i \in w\mathcal{M}^{p_i}_{q_i}(\mathbb{R}^n)$ we have
$$
\biggl\| \prod_{i=1}^{m} f_i \biggr \|_{w\mathcal{M}^{p}_{q}} \leq \prod_{i=1}^{m} \left( \frac{p_i}{p^*} \right)^{\frac{1}{p_i}}
\| f_i \|_{w\mathcal{M}^{p_i}_{q_i}}.
$$
\end{theorem}\vspace{0.2mm}

\noindent
{\sc Proof.} Let $f_i \in w\mathcal{M}^{p_i}_{q_i}$ for $i \in \{1, 2 ,\dots, m\}$, $B = B(a,r) \subseteq \mathbb{R}^n$,
and $\theta>0$. Suppose that $s_0 = \theta$ and $s_m = 1$. Note that for every $s_1, s_2, \dots, s_{m-1}>0$ we have
$$\begin{aligned}[c]
\left| \{ x \in B: \left| f_1 f_2 \cdots f_m(x) \right| > \theta \} \right|
&\leq \sum_{i=1}^{m} \left| \left\{ x \in B: \left| f_i(x) \right| > \frac{s_{i-1}}{s_i} \right\} \right|
\\&\leq \sum_{i=1}^{m} \|f_i\|_{L^{p_i,\infty}(B)}^{p_i}  \left( \frac{s_i}{s_{i-1}} \right)^{p_i}
\\&= \sum_{i=1}^{m} a_i y_i^{p_i}
\end{aligned}$$
with
$$
a_i = \|f_i\|_{L^{p_i,\infty}(B)}^{p_i} \quad \text{and} \quad y_i = \frac{s_i}{s_{i-1}}.
$$
Consequently, for every $y_1, \dots, y_m >0$, we have
$$
\left| \{ x \in B: \left| f_1 f_2 \cdots f_m(x) \right| > \theta \} \right| \leq \sum_{i=1}^m a_i y_i^{p_i}.
$$
By choosing
$$
y_i = \left( m^{-1} \theta^{-p^*} a_i^{-1} \prod_{j=1}^m \left( \frac{p_j}{p^*} \right)^{\frac{p^*}{p_j}}
a_j^{\frac{p^*}{p_j}} \right)^{\frac{1}{p_i}}
$$
we get
$$
\begin{aligned}
\left| \{ x \in B: \left| f_1 f_2 \cdots f_m(x) \right| > \theta \} \right|
&\leq \sum_{i=1}^m a_i y_i^{p_i}
\\&= \sum_{i=1}^m a_i \left( m^{-1} \theta^{-p^*} a_i^{-1} \prod_{j=1}^m \left( \frac{p_j}{p^*}
\right)^{\frac{p^*}{p_j}} a_j^{\frac{p^*}{p_j}} \right)
\\&= \theta^{-p^*} \prod_{j=1}^m \left( \frac{p_j}{p^*} \right)^{\frac{p^*}{p_j}} \|f_j\|_{L^{p_j,\infty}(B)}^{p^*}.
\end{aligned}
$$
Moreover,
$$
\begin{aligned}
\left|B\right|^{\frac{1}{q} - \frac{1}{p*}} \theta \left| \{ x \in B: \left| f_1 f_2 \cdots f_m(x) \right|
> \theta \} \right|^{\frac{1}{p^*}}
&\leq \left|B\right|^{\frac{1}{q} - \frac{1}{p*}} \prod_{i=1}^m \left( \frac{p_i}{p^*} \right)^{\frac{1}{p_i}}
\|f_i\|_{L^{p_i,\infty}(B)}
\\&= \prod_{i=1}^m \left( \frac{p_i}{p^*} \right)^{\frac{1}{p_i}} \left|B\right|^{\frac{1}{q_i} - \frac{1}{p_i}}
\|f_i\|_{L^{p_i,\infty}(B)}
\\&\leq \prod_{i=1}^m \left( \frac{p_i}{p^*} \right)^{\frac{1}{p_i}} \|f_i\|_{w\mathcal{M}^{p_i}_{q_i}}.
\end{aligned}$$
By taking the supremum over all open balls $B$ and $\theta>0$, we have
$$
\biggl\| \prod_{i=1}^{m} f_i \biggr \|_{w\mathcal{M}^{p*}_{q}} \leq \prod_{i=1}^{m}
\left( \frac{p_i}{p^*} \right)^{\frac{1}{p_i}} \| f_i \|_{w\mathcal{M}^{p_i}_{q_i}}.
$$
Hence, by using Theorem \ref{inklusimorreylemah}, we get
$$
\biggl\| \prod_{i=1}^{m} f_i \biggr \|_{w\mathcal{M}^{p}_{q}}
\leq \biggl\| \prod_{i=1}^{m} f_i \biggr \|_{w\mathcal{M}^{p*}_{q}}
\leq \prod_{i=1}^{m} \left( \frac{p_i}{p^*} \right)^{\frac{1}{p_i}} \| f_i \|_{w\mathcal{M}^{p_i}_{q_i}},
$$
which completes the proof.
\qed

\vspace{.5cc}

Next, we prove that the bound in Theorem \ref{holdermorreylemah} is in general sharper than
that in Theorem \ref{holdermorreygunawan}.
To do so, we need the following theorem.

\vspace{.5cc}

\begin{theorem}{[\textbf{AM-GM Weighted Inequality}] \label{amgm}}
Let $m \in \mathbb{N}$. If $x_i, w_i > 0$ for $i=1, 2, \dots, m$, and $w = \sum_{i=1}^m w_i$. Then
$$
\left( \prod_{i=1}^m x_i^{w_i} \right)^{\frac{1}{w}} \leq \sum_{i=1}^m \frac{w_ix_i}{w}.
$$
The equality is attained when $x_i = x_j$ for every $i, j \in \{1,\dots,m\}$.
\end{theorem}\vspace{0.2mm}

\noindent
{\sc Proof.} If $x_i = x_j$ for every $i,j \in \{1, 2, \dots,m\}$, then equality is clearly attained.
Now suppose that $x_i \neq x_j$ for some $i, j \in \{1,\dots,m\}$. Since the natural logarithm is
a strictly concave function, we can use Jensen inequality to get
$$
\ln \left( \prod_{i=1}^m x_i^{w_i} \right)^{\frac{1}{w}} = \ln \left( \prod_{i=1}^m x_i^{\frac{w_i}{w}}
\right) = \sum_{i=1}^m \frac{w_i}{w} \ln \left( x_i \right) < \ln \left(\sum_{i=1}^m \frac{w_ix_i}{w} \right).$$
Hence, we have
$$\left( \prod_{i=1}^m x_i^{w_i} \right)^{\frac{1}{w}} < \sum_{i=1}^m \frac{w_ix_i}{w},$$
as desired.
\qed

\vspace{.5cc}

Now we are ready to prove that our bound is sharper than those in \cite{gunawan2}.

\vspace{.5cc}

\begin{theorem}
Let $m \in \mathbb{N}$, $1 \leq p^* < \infty$, and $1 < p_i < \infty$ for every $i \in \{1,\dots,m\}$.
If $\sum_{i=1}^m \frac{1}{p_i} = \frac{1}{p^*}$, then
$$
\prod_{i=1}^{m} \left( \frac{p_i}{p^*} \right)^{\frac{1}{p_i}} \leq m^{\frac{1}{p^*}} \leq m.
$$
\end{theorem}\vspace{0.2mm}

\noindent
{\sc Proof.} Let $\frac{1}{p^*} = \sum_{i=1}^{m} \frac{1}{p_i}$. By using Theorem \ref{amgm} where
$x_i = \frac{p_i}{p^*}$, $w_i = \frac{1}{p_i}$, and $w = \frac{1}{p^*}$, we get
$$
\left( \prod_{i=1}^{m} \left( \frac{p_i}{p^*} \right)^{\frac{1}{p_i}} \right)^{p^*} \leq
\sum_{i=1}^m \frac{1}{p_i} \frac{p_i}{p^*} p^* = m.
$$
Hence,
$$
\prod_{i=1}^{m} \left( \frac{p_i}{p^*} \right)^{\frac{1}{p_i}} \leq m^{\frac{1}{p^*}} \leq m,
$$
as stated.
\qed

\vspace{1.5cc}

\section{CONCLUDING REMARKS}

In this note we already proved that generalized H\"{o}lder's inequality in weak Morrey spaces with
bound $\prod_{i=1}^{m} \left( \frac{p_i}{p^*} \right)^{\frac{1}{p_i}}$ and this bound is in general
sharper than the one obtained in \cite{gunawan2}.

However, we still do not know whether the bound is really sharp, that is, we still do not know whether
there is a function $f_i \in w\mathcal{M}^{p_i}_{q_i}(\mathbb{R}^n)$, $i=1, \dots, m$ for some $m \geq 2$,
such that
$$
\biggl\| \prod_{i=1}^{m} f_i \biggr \|_{w\mathcal{M}^{p}_{q}} = \prod_{i=1}^{m} \left( \frac{p_i}{p^*}
\right)^{\frac{1}{p_i}} \| f_i \|_{w\mathcal{M}^{p_i}_{q_i}}.
$$

\vspace{1.5cc}

\noindent{\bf Acknowledgement.} The results in this note have been presented at Mid Year School on
Analysis, Geometry, and Applications (MYSAGA) 2018. We would like to thank Ms.~Ifronika for her comments
on the earlier version of this note. 


\end{document}